\documentclass[11pt]{amsart}
\usepackage[margin=2cm]{geometry}  
\newtheorem{theorem}{Theorem}                
\usepackage{amssymb}
\usepackage{amsmath}
\usepackage{xcolor}
\usepackage{url}

\usepackage{listings}
\lstdefinelanguage{Sage}[]{Python}
{morekeywords={False,sage,True},sensitive=true}
\definecolor{dblackcolor}{rgb}{0.0,0.0,0.0}
\definecolor{dbluecolor}{rgb}{0.01,0.02,0.7}
\definecolor{dgreencolor}{rgb}{0.2,0.4,0.0}
\definecolor{dgraycolor}{rgb}{0.30,0.3,0.30}

\begin{document}

\title{Unitary and Modular Fourier Transforms of the Symmetric Group over Finite Fields}
\author[J.~Walters]{Jackson Walters}
\address[J.~Walters]{Arlington, VA}
\email{jacksonwalters@gmail.com}
\urladdr{https://jacksonwalters.com}


\begin{abstract}
We give explicit constructions of Fourier transforms for cyclic and symmetric groups over finite fields. For a cyclic group in the modular case $p\mid N$, we factor $x^N-1$ and use B\'{e}zout's identity to construct a Fourier transform and its inverse. For the symmetric group in the modular case $p\mid n!$, we compute the Peirce decomposition from central primitive orthogonal idempotents, obtaining a change-of-basis matrix which generalizes the DFT. We also construct unitary realizations of the Specht representations of the symmetric group over finite fields: we compute an invariant Hermitian form and factor its Gram matrix by an extended Cholesky decomposition. These representations supply the matrix coefficients for a unitary DFT in the semisimple case. All constructions are implemented in SageMath.
\end{abstract}

\maketitle

\section{Introduction}

The Fourier transform is an old and powerful tool used in variety of applications, from signal processing to fundamental physics. It is central in mathematics. At its core, it is a change of basis that allows us to view a signal changing in time in terms of static frequencies.  This article aims to build intuition for an extension of the usual Fourier transform that takes place over the field of real or complex numbers to a more general field of potentially positive characteristic. Further, we implement the so called ``modular discrete Fourier transform" in SageMath \cite{sage-math} for the cyclic group and for the symmetric group $S_n$. \\

The paper has three main constructions. First, we give a concrete DFT for cyclic groups over finite fields, including the modular case. Second, we use the block decomposition of the modular group algebra to define and compute a DFT for the symmetric group. Third, we construct unitary matrix realizations of the symmetric group over finite fields by computing invariant Hermitian forms and their extended Cholesky decompositions. \\

The modular discrete Fourier transform for a field $K$ with $\text{char}(K) > 0$ uses the isomorphism from the group algebra $K\left[G\right]$ to its block decomposition $\bigoplus_i K\left[G\right]e_i$, where $1=\sum_i e_i$ is a decomposition of the identity into primitive orthogonal idempotents. This generalizes the usual Fourier transform which uses Maschke's theorem to decompose $K\left[G\right] = \bigoplus_i \text{End}(V_i)$ for irreducible representations $V_i$ of $G$ so long as $\text{char}(K) \nmid |G|$. \\

The article will be primarily concerned with the Fourier transform for finite groups, the simplest example being abelian groups which gives rise to the discrete Fourier transform (DFT). This is a map from $\mathbb{C}^N \rightarrow \mathbb{C}^N$ which can be represented by an $N\times N$ matrix. The DFT has wide applicability in engineering, and is amenable to implementation on a computer. \\

The DFT generalizes readily to the realm of quantum computing. We view a collection of $n$ qubits $c_0 | 0 \rangle + c_1 | 1 \rangle \in \mathbb{C}^2$ which are concatenated in space to correspond to the Hilbert space $\mathcal{H} \cong \left(\mathbb{C}| 0 \rangle + \mathbb{C}| 1 \rangle \right)^{\otimes n} \cong \mathbb{C}^{2^n}$. Then the Fourier transform is a unitary matrix $U: \mathcal{H} \rightarrow \mathcal{H}$, where we label the eigenstates using binary notation, e.g. for $n=3$, $| 0 \rangle \otimes | 1 \rangle \otimes | 1 \rangle = |011\rangle$ can be written as $| 3 \rangle$ in the familiar decimal, and we interpret the label modulo $N=2^n$. We can then define maps such as $| x \rangle \mapsto | a^x (\text{mod} \, N) \rangle$ which is useful in applications like Shor's algorithm \cite{peter-shor-94}, where we use the inverse quantum Fourier transform to extract the periodicity of that modular exponentiation viewed as a function $\mathbb{Z}/N\mathbb{Z} \rightarrow \mathbb{Z}/N\mathbb{Z}$. This technique is referred to as ``Fourier sampling." \\

\section*{Acknowledgements}

I would like to thank Ezra Brown, Noam Elkies, Dmitrii Pasechnik, Robert Pollack, Sridhar Ramesh, Steve Rosenberg, Travis Scrimshaw, and Jared Weinstein for helpful comments and conversation.

\section{The Discrete Fourier Transform}

It is important that the base field, $\mathbb{C}$, is of characteristic zero. This means that we can always invert the order of the cyclic group, $N$, to define the inverse Fourier transform. However, we have thus far suppressed the representation theory and harmonic analysis. Note that $\mathbb{Z}/N\mathbb{Z}$ is an abelian group, so its representations are one dimensional. So we are looking at maps $\chi: \mathbb{Z}/N\mathbb{Z} \rightarrow \mathbb{C}$. This means that $\chi(N) = \chi(1)^N =1$, which means $\chi(1)$ is an $N^{th}$ root of unity. Since these are one dimensional, they are irreducible. Thus, there are $N$ representations given by $\chi_k(n)=\exp(\frac{-2\pi i k n}{N})$. \\

Note that $\chi_k$ are functions on the group $\mathbb{Z}/N\mathbb{Z}$. Thus, for each group element $n \in \{0,\ldots,N-1\}$ we have a complex number $\chi_k(n) \in \mathbb{C}$. Recall that the space of class functions on a group form a vector space. Since the group is abelian, each conjugacy class has one element, so $\chi_k$ may be viewed as a vector in $V^* \cong \mathbb{C}^N$. Separately, we have the group alegbra $V = \mathbb{C}\left[\mathbb{Z}/N\mathbb{Z}\right] \cong \mathbb{C}^N$ as a vector space. We can multiply these ``vectors" by viewing group elements as a basis and extending the multiplication as we do with polynomials. Label the irreducible representations $V_k \cong \mathbb{C}$. Each element $n \in \mathbb{Z}/N\mathbb{Z}$ may be viewed as an endomorphism $n: V \rightarrow V$ by multiplication. That is, if $s=\sum_j c_j j \in V$ then $n\cdot s = \sum_k c_j (n \cdot j)$ where the product $n\cdot j = n+j \mod N$ is in $\mathbb{Z}/N\mathbb{Z}$. Therefore we have a linear action of $\mathbb{Z}/N\mathbb{Z}$ on $V$, which is invertible since $-n$ is the multiplicative inverse of $n$. Thus, each $n \in GL(V)$, so this is a representation of the group. $V$ can be decomposed into irreducible representations $V \cong \bigoplus_k V_k$ and the action of $n$ diagonalized as $n_k \in End(V_k)$. We thus have the decomposition $V=\mathbb{C}\left[ \mathbb{Z}/N\mathbb{Z}\right] \cong \bigoplus_k \text{End}(V_k) \cong \bigoplus_k \mathbb{C} \chi_k = V^*$. \\

So how does the Fourier transform look in this interpretation? Consider the element $x = \left(x_1, \ldots, x_N\right) \in \mathbb{C}^N$. We map it to $\hat{x}_k = \sum_{n=0}^{N-1} x_n \exp\left(\frac{-2\pi i k n}{N}\right)$. \\

Let us rewrite this formula a bit. 

$$\hat{x}_k = \sum_{n=0}^{N-1} x_n \exp\left(\frac{-2\pi k}{N}\right)^n = \sum_{n=0}^{N-1} x_n \chi_k(n)$$

We see that our ``Fourier coefficients" are just sums of the components of $x$ against the values of the irreducible characters. We obtain $N$ coefficients, and this is exactly the discrete Fourier transform for the group $\mathbb{Z}/N\mathbb{Z}$. This formula has the nice property of being invertible, with inverse given by $x_n = \frac{1}{N}\sum_{k=0}^{N-1}\hat{x}_k \chi_k(-n)$, whose proof uses the orthogonality of characters. \\

This is a bit of a special case, in which the number of representations equals the order of the group, i.e. every representation is degree 1. This only occurs for abelian groups in algebraically closed fields of characteristic 0. For nonabelian groups, we need to use matrices. \\

\section{Fourier Transform for Non-Abelian Groups} \label{fourier-transform-nonabelian-groups}

Let $G$ be a finite group and let $\widehat{G}$ denote the set of irreducible representations of $G$, also a finite set. Over a field of characteristic zero, we have Maschke's theorem:

\begin{theorem} 
Let $G$ be a finite group and $K$ a field whose characteristic does not divide the order of $G$. Then $K\left[G\right]$, the group algebra of $G$, is semisimple.
\end{theorem}

This theorem tells us that as a $K\left[G\right]$-module we may decompose $K\left[G\right]$ into a direct sum of simple $K\left[G\right]$-modules. That is, $K\left[G\right] = \bigoplus_i M_i$. It turns out these simple $M_i$ are just the endomorphism algebras of irreducible representations, i.e. $K\left[G\right] = \bigoplus_i End(V_i)$. Thus, an element of the group algebra may be viewed as a ``block diagonal" matrix acting on itself, with each block corresponding to an irreducible representation. This is a special case of Wedderburn's theorem. \\

By analogy, we may form the ``matrix Fourier coefficients" as sums of coefficients times irreducible representations. Let $f: G \rightarrow K$ be a function on the group, or equivalently an element of the group algebra $f \in K\left[G\right]$. Define for $\rho \in \hat{G}$

$$\hat{f}(\rho) = \sum_{a \in G}f(a)\rho(a)$$

\noindent where $f(a) \in K$ and each $\rho(a)$ is a $d_\rho \times d_\rho$ matrix. This formula is invertible yielding

$$f(a) = \frac{1}{|G|}\sum_{g \in \hat{G}} d_\rho \text{Tr}(\rho(a^{-1})\hat{f}(\rho)).$$

It is clear we require $\text{char}(K) \nmid |G|$, both in order to comply with the hypothesis of Mashke's Theorem, and to avoid division by zero. We need to know the dimension numbers of irreducible representations, $d_\rho$. We are dealing with matrices, so we are taking the matrix trace now in the sum. \\

\section{The Symmetric Group over $\mathbb{C}$}

For the symmetric group over $\mathbb{C}$, Maschke's theorem applies. The irreducible representations are labeled by Young tableaux $\lambda$ corresponding to partitions of $m$. A \emph{tabloid} is an equivalence class of Young tableaux where two labelings are equivalent if one is obtained from the other by permuting the entries of each row. For each Young tableau $T$ of shape $\lambda$ let $\{T\}$ be the corresponding tabloid. For each $\lambda$, we obtain a simple modules $V^\lambda$ with basis given by

$$E_T = \sum_{\sigma \in Q_T} \epsilon(\sigma)\{\sigma(T)\} \in V^\lambda$$

\noindent where $Q_T$ is the subgroup of permutations, preserving (as sets) all columns of $T$ and $\epsilon(\sigma)$ is the sign of the permutation $\sigma$. The Specht module of the partition $\lambda$ is the module generated by the elements $E_T$ as $T$ runs through all tableaux of shape $\lambda$. The dimension of the Specht module $V^\lambda$ is the number of standard Young tableaux of shape $\lambda$. It is given by the hook length formula. Over fields of characteristic 0 the Specht modules are irreducible and form a complete set of irreducible representations. \\

We now have everything we need to compute the non-abelian Fourier transform, which is implemented in \cite{dft-finite-groups} for the symmetric group (though it generalizes to any finite group given a labeling for the representations).

\begin{lstlisting}
#define a permutation group of size n!
n=4
G=SymmetricGroup(n)

#define the Fourier transform at the representation spc
#which is the Specht module corresponding to partition
def f_hat(f,partition):
    spc = SymmetricGroupRepresentation(partition, 'specht')
    return sum(f(sigma)*spc.representation_matrix(Permutation(sigma)) for sigma in G)

def fourier_coeffs(f):
    return {partition:f_hat(f,partition) for partition in Partitions(n)}
    
#define the inverse Fourier transform
#given a collection of matrices for each representation rho_i
def fourier_inv(sigma,f_hat):
    value=0
    for partition in Partitions(n):
        spc = SymmetricGroupRepresentation(partition, 'specht')
        repn_mat = spc.representation_matrix(Permutation(sigma).inverse())
        d_part = repn_mat.ncols() #should extract directly from parition of rep'n
        fourier_coeff = f_hat[partition] #fourier coefficient
        value += d_part*(repn_mat*fourier_coeff).trace()
    return (1/G.order())*value

\end{lstlisting}

\section{When Maschke's Theorem Fails, $p \mid |G|$}

There is a well-known isomorphism of categories between $K\left[G\right]$-modules are representations of the group $G$ over the field $K$. The ``regular representation" where $k\left[G\right]$ acts on itself by multiplication decomposes into a sum of irreducible representations. When the characteristic of the field $p$ does not divide the order of the group $|G|$ we can use Maschke's theorem as stated above. The key element of the proof is the construction of an ``average" or invariant inner product, $\langle x,y\rangle_G = \frac{1}{|G|}\sum_{g \in G} \langle g x,g y\rangle$, where juxtaposition denotes the group action in a representation. Note that typically we require $K$ be algebraically closed, since otherwise a more complicated statement holds: $K\left[G\right]$ is a product of matrix algebras over division rings over $K$. In the case of finite fields, any division ring over $F_p$ must be a field by Wedderburn's little theorem, and the only finite fields of characteristic $p$ are the extensions $F_q$ where $q=p^\nu$. \\

When $p \mid |G|$, the situation is more complicated and is the subject of modular representation theory. In our case, we are primarily concerned with cyclic groups and the symmetric group $S_m$. In general, we no longer have Maschke's theorem giving us a decomposition of the group algebra into simple modules. Instead, we have a decomposition into ``blocks",

$$K\left[G\right] = \bigoplus_i K\left[G\right]e_i$$

\noindent where $\sum_i e_i = 1$ is a decomposition of the identity into central primitive orthogonal idempotents. In certain cases these idempotents are readily computable.

\section{The Discrete Fourier Transform over a Ring}

We implement this procedure in SageMath code. \cite{sage-math} \cite{dft-finite-groups}  \\

The discrete Fourier transform as developed above works over the complex numbers. However, the basic formula only uses a homomorphism from a cyclic group to a field, $\chi_k(n) = \exp(\frac{2 \pi n}{N})$. We just need an element $\alpha$ such that $\alpha^N=1$, and $\sum_{j=0}^{N-1} \alpha^{jk} = 0$ for $0 \le k < n$. Then we may form a map $\left(v_0, \ldots, v_{N-1}\right) \mapsto \left(f_0, \ldots, f_{N-1}\right)$ where

$$f_k = \sum_{j=0}^{N-1} v_j \alpha^{j k}$$

If $R$ is an integral domain (which includes fields), then we only require $\alpha$ be a primitive $N^{th}$ root of unity, i.e. $\alpha^k \ne 1$ for $1 \le k < N$. \\

When $R = \mathbb{F}_p$, the multiplicative group $\mathbb{F}_p^{\times} \cong \mathbb{Z}/(N-1)\mathbb{Z}$ is cyclic. By Fermat's little theorem, $\alpha^{p-1} = 1$. Thus, we require $N | (p-1)$. This implies that $p > N+1$, which means that the characteristic we choose must be bigger than the size of the list we are transforming. This is quite a constraint! \\

We also require that $p \nmid N$, to ensure that $N$ is invertible in $\mathbb{F}_p$. Otherwise we cannot define the inverse Fourier transform. \\

\begin{lstlisting}
n=16; p=17
#finite field of size p
K = GF(p) 
assert K(n) != K(0) #ensure n is invertible
assert n.divides(p-1) #ensure a primitive n-th root of unity exists
            
#list to be transformed
v = [K(i) for i in range(n)]
alpha = K.zeta(n)

#DFT over a ring with primitive root alpha
def fourier_transform(v,alpha):
    return [sum(v[j]*alpha**(j*k) for j in range(len(v))) for k in range(len(v))]
f=fourier_transform(v,alpha); f
[1, 8, 2, 15, 7, 4, 6, 5, 9, 13, 12, 14, 11, 3, 16, 10]
#define the inverse Fourier transform
def inverse_fourier_transform(f,alpha):
    return [K(1/len(f))*sum(f[k]*alpha**(-j*k) \ 
    	      for k in range(len(f))) for j in range(len(f))]
inverse_fourier_transform(f,alpha)
[0, 1, 2, 3, 4, 5, 6, 7, 8, 9, 10, 11, 12, 13, 14, 15]
\end{lstlisting}

\section{Cyclic Groups over $\mathbb{F}_p$} \label{cyclic-groups-finite-field}

We now consider the group algebra $\mathbb{F}_p \left[\mathbb{Z}/N\mathbb{Z}\right]$. \\

The question of how to classify all finite dimensional representations of $\mathbb{Z}/N\mathbb{Z}$ over an arbitrary field $F$ can be studied using the structure theorem for finitely-generated modules over a principal ideal domain, in this case $F\left[x\right]$. The structure theorem asserts that any finitely-generated module is uniquely a finite direct sum of modules of the form $F\left[x\right]/p(x)^r$ where $p \in F\left[x\right]$ is irreducible and $r$ is a non-negative integer. \\

If $T$ is an operator acting on $F^k$ for some $n$, then $F^k$ becomes a finitely-generated module over $F\left[x\right]$ with $x$ acting by $T$. This gives a representation of the cyclic group $\mathbb{Z}/N\mathbb{Z}$ if and only if $T^N=1$, in which case the summands $F\left[x\right]/q(x)^r$ in the decomposition of $F^k$ must have the property that $q(x)^r | (x^N -1)$. \\

\subsection{Case $p \nmid N$} If $F$ has characteristic 0 or has characteristic $p$ and $p \nmid N$, then $x^N - 1$ is separable over $F$, hence $r \le 1$ and $F^k$ is a direct sum of irreducible representations, all of which are of the form $F\left[T\right]/q(T)$ where $q$ is an irreducible factor of $x^N-1$ over $F$. \\

In module language, that is $F\left[\mathbb{Z}/N\mathbb{Z}\right] \cong F\left[x\right]/(x^N - 1) \cong \bigoplus_k F\left[x\right]/q_k(x)$, where $q_k(x)$ are irreducible polynomials over $F$. \\

What we are really doing is factoring $x^N-1$ in $F_p\left[x\right]$. We can first factor this into a product of cyclotomic polynomials $\Phi_d(x)$ in $\mathbb{Z}\left[x\right]$. That is, $x^N-1 = \prod_{d|N} \Phi_d(x)$. Further, we wish to factor $\Phi_d(x)$ in $F_p$. This is equivalent to factoring the the prime ideal $(p)$ in $\mathbb{Z}\left[\zeta\right] = \mathbb{Z}\left[x\right]/\Phi_d(x)$. It turns out $\Phi_d = P_1 \ldots P_g$ where $P_i$ each have residue degree $f$, and $fg = \varphi(d)$ where $\varphi$ is Euler's totient function, and $f$ is the order of $p \mod d$. \\

Note that solutions to $x^N-1 = 0$ in $F_p$ are just roots of unity in a finite field. Since $x^{p-1} -1 = 0$ in any field of order $p$, we require that $N|(p-1)$ to have solutions. Further, when there are solutions there are $gcd(N,p-1)$ of them. 

\subsubsection{Examples} 

$$N=6, p=7$$

In this case, $N|p-1$, so we expect to have a primitive $6^{th}$ root of unity in $F_7$. Once we have such a root $\alpha$, all $k^{th}$ powers of $\alpha$ are roots of $x^N-1$ as well, so we factor 

$$x^6-1 = (x+1)(x+2)(x+3)(x+4)(x+5)(x+6)$$

$$N = 10, p=7$$

$$x^{10} - 1 = \prod_{d|10}\Phi_d(x) = \prod_{1,2,5,10}\Phi_d(x) = (x-1)(x+1)(x^4+x^3+x^2+x+1)(x^4-x^3+x^2-x+1)$$

How do $\Phi_5(x)$ and $\Phi_{10}(x)$ factor in $F_p\left[x\right]$? Well, $\varphi(5)=5-1=4$. $7 \equiv 2 \mod 5$, so $7$ has order 4 $\mod 5$. Thus, when $d=5$, $f=4$, so $g=1$, thus there is only one polynomial in the factorization, the cyclotomic polynomial. When $d=10$, $\varphi(10)=4$. $7^2 = 49 \equiv -1 \mod 10$, so $7^4 \equiv 1 \mod 10$, thus 7 has order 4 $\mod 10$, so $f=4$, and $g=1$. Therefore our factorization is complete. \\

Note that we may use the Chinese remainder theorem to obtain an isomorphism to the decomposition into summands: \cite{dft-finite-groups}

$$F_p\left[x\right] / (x^N-1) \cong \bigoplus_{d|N,i} F_p\left[x\right] / (P_i(x))$$

\noindent The map is just given by $f \mapsto f \mod P_i(x)$ for each of the $g$ $P_i$'s occuring in the factorization of the cyclotomic polynomials. The inverse map is given by multiplying the factors in the variable $x$, and taking the residue modulo $x^N-1$. \\

We may want to find an extension $F_q$ of $F_p$ where $q=p^\nu$ that splits $x^N-1$. This is always possible, since we just need to find a primitive $N^{th}$ root, which occurs when $N|(q-1)$. In this case, we have

$$x^N-1 = (x-\alpha^{0})(x-\alpha^1)\ldots(x-\alpha^{N-1})$$

\noindent where $\alpha^{jN} = 1$ for all $0 \le j \le N-1$. In this case, the Fourier transform becomes the usual map 

$$f_k: x \mapsto \sum_j x_j \alpha^{jk}$$

\begin{lstlisting}

#h is an element of F_p[C_N]
#that is, h = h0+h1*x+h2x^2+...+h_{N-1}x^N-1
#we allow h as a list of N numbers modulo p
#h = [h0,h1,h2,...,h_{N-1}]
def discrete_fourier_transform(h,p,splitting_field=False):
    #length of list is size of N
    N = len(h)
    #define the polynomial ring F_p[x]
    R = PolynomialRing(GF(p),'x')
    #name the generator x an element of R
    x = R.0
    #define the polynomial x^N-1
    f = x**N-1; assert f in R
    if splitting_field:
        K.<a> = f.splitting_field()
        #define the polynomial ring over extended base field
        R = PolynomialRing(K,'x')
        #name the generator x an element of R
        x = R.0
        #define the polynomial x^N-1
        f = x**N-1; assert f in R
    #define the quotient ring F_p[x]/(x^N-1)
    S = R.quotient(x^N - 1, 'z')
    #transform the list of coefficients of h into a polynomial in R=F_p[x]
    h = sum(h[i]*x**i for i in range(N)); assert h in S
    #factor f in F_p[x], save as list of factors and multiplicities
    f_factors = list(f.factor())
    #implement the Chinese remainder theorem mapping 
    #S=F_p[x]/(x^N-1) --> \prod_i R/(factor_i^mult_i)
    h_transform=[list(R.quotient(f_factors[i][0]**f_factors[i][1])(h)) \
    	for i in range(len(f_factors))]
    return h_transform
    
def inv_discrete_fourier_transform(hhat,p,splitting_field=False):
    N = sum(len(l) for l in hhat)
    #define the polynomial ring F_p[x]
    R = PolynomialRing(GF(p),'x')
    #name the generator x an element of R
    x = R.0
    #define the polynomial x^N-1
    f = x**N-1; assert f in R
    if splitting_field:
        K.<a> = f.splitting_field()
        #define the polynomial ring over extended base field
        R = PolynomialRing(K,'x')
        #name the generator x an element of R
        x = R.0
        #define the polynomial x^N-1
        f = x**N-1; assert f in R
    S = R.quotient(x^N - 1, 'x')
    f_factors = list(f.factor())
    #perform inverse of Chinese remainder theorem
    #for each modulus N_i = N/n_i, where n_i is the modulus of each factor
    #Bezout's theorem applies, so we get M_i*N_i + m_i*n_i = 1
    #a solution x = \sum_{i=1}^k a_i*M_i*N_i, where a_i are the remainders
    n = [f_factors[i][0]**f_factors[i][1] for i in range(len(f_factors))]
    #get coefficients M_i, m_i from N_i, n_i
    M = [xgcd(f/n[i],n[i])[1] for i in range(len(n))]
    #get remainders as polynomials in R
    a = [sum(hhat[i][j]*x**j for j in range(len(hhat[i]))) for i in range(len(hhat))]
    inv_transform = sum(a[i]*M[i]*(f/n[i]) for i in range(len(a)))
    return list(S(inv_transform))

\end{lstlisting}

\vspace{5mm}

\subsection{Case $p \mid N$} Writing $N = p^s \ell$ where $p \nmid \ell$ we have

$$x^N - 1 = (x^\ell - 1)^{p^s}$$

\noindent so it follows that $r \le p^s$ (but now it is possible that $r>1$). If $r>1$, then the corresponding representation $F\left[T\right]/q(T)^r$ is indecomposable and not irreducible, where $q$ is an irreducible factor of $x^\ell - 1$ over $F$. The irreducible representations occur precisely when $r=1$. \\

If $V$ is an irreducible representation of $\mathbb{Z}/N\mathbb{Z}$, and $T: V \rightarrow V$ is the action of a generator, then 

$$T^{p^s \ell} - 1 = (T^\ell - 1)^{p^s} = 0$$

Thus $T^\ell - 1$ is an intertwining operator which is not invertible, so by Schur's lemma it is equal to zero. \\

In this case, we may still factor $x^\ell - 1$ as before, we just obtain multiplicities on all factors, so the representations are now indecomposable rather than irreducible. Thus, we have the isomorphism

$$F_p\left[x\right] / (x^N-1) \cong F_p\left[x\right] / (P_i(x)^{p^s}).$$ 

\section{The Symmetric Group over $\mathbb{F}_p$}

The symmetric group algebra can be decomposed into blocks $F_p\left[S_N\right] = \bigoplus_i F_p\left[S_N\right] e_i$ where $1=\sum_i e_i$ is a decomposition of the identity into primitive, central, orthogonal idempotents. This is known as the Peirce decomposition. \\

Though we do not make use of them, the simple modules $D^\lambda = S^\lambda / \left(S^\lambda \cap (S^\lambda)^{\perp}\right)$ occur for each $p$-regular partition $\lambda$, and form a complete traversal. They are computed in SageMath 10.4+ \cite{sage-math} for small $n$, and much more extensive tables are found in \cite[pg. 98]{james-book-78}; see also \cite{james-paper-88}. \\

The idempotents are computed by Murphy \cite{murphy-idem-83} as sums involving the Jucys--Murphy elements, and equivalence classes of partitions based on their $p$-core. The $p$-core of a partition is obtained by removing as many rim $p$-hooks from the edge of a partition as possible. The number removed is the weight. $p$-cores $\gamma$ parameterize the blocks of the symmetric group. \cite{wildon-notes} \\

We include an explicit computation of the idempotents in \cite{dft-finite-groups}, and they can be computed much more quickly using Sage via 

\begin{lstlisting}
p=3; n=4;
SGA = SymmetricGroupAlgebra(GF(p),n)
print(SGA.central_orthogonal_idempotents())
\end{lstlisting}

Once we have computed the idempotents, we can write $F_p\left[S_N\right]e_i=\text{span}_{F_p}\{\sigma e_i | \sigma \in S_N\}$ to find a basis for the block. Putting these bases together, we obtain a second basis for the symmetric group algebra different from the basis $\{\sigma | \sigma \in S_N\}$. The change of basis matrix is the modular Fourier transform \cite{dft-finite-groups}. \\

\begin{lstlisting}
#implements modular Fourier transform
#project v onto each block U_i = F_p[S_n]*e_i 
#using \pi_i: v |--> v*e_i as a projection
#this is just a change of basis
def modular_fourier_transform(p,n):
    #instantiate group algebra
    SGA_GFp_n = SymmetricGroupAlgebra(GF(p),n)
    #compute the primitive central orthogonal idempotents
    idempotents = SGA_GFp_n.central_orthogonal_idempotents()
    #create a spanning set for the block corresponding to an idempotent
    spanning_set = lambda idem: [SGA_GFp_n(sigma)*idem for sigma in SGA_GFp_n.group()]
    #compute the blocks as submodules of the given spanning set
    blocks = [SGA_GFp_n.submodule(spanning_set(idem)) for idem in idempotents]
    #compute the list of basis vectors lifed to the SGA from each block
    block_decomposition_basis = flatten([[u.lift() \
    	for u in block.basis()] for block in blocks])
    #the elements of the symmetric group are ordered
    #giving the map from the standard basis
    sym_group_list = list(SGA_GFp_n.group())
    change_of_basis_matrix = []
    for b in block_decomposition_basis:
        coord_vector = [0]*len(sym_group_list)
        for pair in list(b):
            coord_vector[sym_group_list.index(pair[0])] = pair[1]
        change_of_basis_matrix.append(coord_vector)
    return matrix(GF(p),change_of_basis_matrix).transpose()
\end{lstlisting}

\subsection{Examples}

\subsubsection{$p=2$, $n=3$}

\begin{lstlisting}
modular_fourier_transform(2,3)
[1 0 0 0 1 0]
[0 1 0 0 0 1]
[0 0 1 0 0 1]
[0 0 0 1 1 0]
[1 0 0 1 1 0]
[0 1 1 0 0 1]
\end{lstlisting}

\subsubsection{$p=3$, $n=4$}

\begin{lstlisting}
print(modular_fourier_transform(3,4).str())
[1 0 0 0 0 0 0 0 0 1 0 0 0 0 0 0 0 0 1 0 0 0 0 0]
[0 1 0 0 0 0 0 0 0 0 1 0 0 0 0 0 0 0 0 1 0 0 0 0]
[0 0 1 0 0 0 0 0 0 0 0 1 0 0 0 0 0 0 0 0 1 0 0 0]
[0 0 0 1 0 0 0 0 0 0 0 0 1 0 0 0 0 0 0 0 0 1 0 0]
[0 0 0 0 1 0 0 0 0 0 0 0 0 1 0 0 0 0 0 0 0 0 1 0]
[1 2 2 1 1 0 0 0 0 2 2 2 2 2 0 0 0 0 0 0 0 0 0 1]
[0 0 0 0 0 1 0 0 0 0 0 0 0 0 1 0 0 0 0 1 0 0 0 0]
[2 1 0 0 0 1 0 0 0 2 2 0 0 0 2 0 0 0 1 0 0 0 0 0]
[0 0 0 0 0 0 1 0 0 0 0 0 0 0 0 1 0 0 0 0 0 0 1 0]
[0 0 0 0 0 0 0 1 0 0 0 0 0 0 0 0 1 0 0 0 0 0 0 1]
[0 0 2 0 1 0 1 0 0 0 0 2 0 2 0 2 0 0 0 0 1 0 0 0]
[1 2 2 0 1 0 0 1 0 1 1 1 0 1 0 0 2 0 0 0 0 1 0 0]
[0 0 0 0 0 0 0 0 1 0 0 0 0 0 0 0 0 1 0 0 0 1 0 0]
[0 0 2 1 0 0 0 0 1 0 0 2 2 0 0 0 0 2 0 0 1 0 0 0]
[1 0 2 0 0 2 1 0 1 2 0 2 0 0 2 2 0 2 0 0 0 0 0 1]
[1 0 2 0 0 2 0 1 1 1 0 1 0 0 1 0 2 1 0 0 0 0 1 0]
[1 2 1 1 1 2 1 0 1 1 1 2 1 1 1 1 0 1 1 0 0 0 0 0]
[1 2 2 0 1 2 0 1 1 2 2 2 0 2 2 0 1 2 0 1 0 0 0 0]
[1 1 2 2 2 1 2 2 2 2 1 2 1 1 1 1 2 1 0 0 0 0 0 1]
[2 0 1 0 2 1 2 2 2 2 0 2 0 2 2 2 1 2 0 0 0 0 1 0]
[2 1 1 2 2 0 0 2 2 2 2 2 2 2 0 0 1 2 0 0 0 1 0 0]
[0 0 1 2 2 0 2 0 2 0 0 1 1 1 0 1 0 1 0 0 1 0 0 0]
[2 0 1 0 2 0 0 2 2 1 0 1 0 1 0 0 2 1 0 1 0 0 0 0]
[2 0 2 2 2 0 2 0 2 2 0 1 2 2 0 2 0 2 1 0 0 0 0 0]
\end{lstlisting}

\subsection{$p=7$, $n=4$} Here $p \nmid n!$, and we obtain a matrix of (almost) entirely 0's and $\pm 1$'s. 

\begin{lstlisting}
#note: this matrix has (almost) only has 1's and -1's
print(modular_fourier_transform(7,4).str())
[1 1 0 0 0 0 0 0 0 0 1 0 0 0 1 0 0 0 0 0 0 0 0 1]
[1 0 1 0 0 0 0 0 0 0 0 1 0 0 0 1 0 0 0 0 0 0 0 6]
[1 0 0 1 0 0 0 0 0 0 0 0 1 0 0 0 1 0 0 0 0 0 0 6]
[1 0 0 0 1 0 0 0 0 0 0 0 0 1 0 0 0 1 0 0 0 0 0 1]
[1 0 0 0 0 1 0 0 0 0 6 0 0 6 0 0 0 0 1 0 0 0 0 1]
[1 1 6 6 1 1 0 0 0 0 0 6 6 0 6 6 6 6 6 0 0 0 0 6]
[1 0 0 0 0 0 1 0 0 0 0 1 0 0 0 0 0 0 0 1 0 0 0 6]
[1 6 1 0 0 0 1 0 0 0 1 0 0 0 6 6 0 0 0 6 0 0 0 1]
[1 0 0 0 0 0 0 1 0 0 6 0 0 6 0 0 0 0 0 0 1 0 0 1]
[1 0 0 0 0 0 0 0 1 0 0 6 6 0 0 0 0 0 0 0 0 1 0 6]
[1 0 0 6 0 1 0 1 0 0 0 0 1 0 0 0 6 0 6 0 6 0 0 6]
[1 1 6 6 0 1 0 0 1 0 0 0 0 1 1 1 1 0 1 0 0 6 0 1]
[1 0 0 0 0 0 0 0 0 1 0 0 0 1 0 0 0 0 0 0 0 0 1 1]
[1 0 0 6 1 0 0 0 0 1 0 0 1 0 0 0 6 6 0 0 0 0 6 6]
[1 1 0 6 0 0 6 1 0 1 0 6 6 0 6 0 6 0 0 6 6 0 6 6]
[1 1 0 6 0 0 6 0 1 1 6 0 0 6 1 0 1 0 0 1 0 6 1 1]
[1 1 6 5 1 1 6 1 0 1 1 0 0 0 1 1 2 1 1 1 1 0 1 1]
[1 1 6 6 0 1 6 0 1 1 0 1 0 0 6 6 6 0 6 6 0 1 6 6]
[1 5 1 2 6 6 1 6 6 6 0 6 6 0 2 1 2 1 1 1 1 6 1 6]
[1 6 0 1 0 6 1 6 6 6 6 0 0 6 6 0 6 0 6 6 6 1 6 1]
[1 6 1 1 6 6 0 0 6 6 0 0 0 1 6 6 6 6 6 0 0 1 6 1]
[1 0 0 1 6 6 0 6 0 6 0 0 1 0 0 0 1 1 1 0 1 0 1 6]
[1 6 0 1 0 6 0 0 6 6 0 1 0 0 1 0 1 0 1 0 0 6 1 6]
[1 6 0 2 6 6 0 6 0 6 1 0 0 0 6 0 5 6 6 0 6 0 6 1]
\end{lstlisting}

\section{Unitarity for the cyclic group}

The discrete Fourier transform for the cyclic group over the complex numbers is not unitary when given by $A = \left(\omega^{ij} \right)_{0 \le i, j < N}$ where $\omega = \exp{2\pi i / N}$. Instead, $AA^*=NI_N$. Multiplying by $\frac{1}{\sqrt{N}}$ gives the unitary matrix $U = \frac{1}{\sqrt{N}}A$. \\

Over finite fields, the situation is slightly more delicate but is in complete analogy to the case over the complex numbers. First we need a notion of conjugation appropriate for finite fields. Suppose we are working over a finite field $F_q$. If $q$ is square, then $q=(q')^2$ and if not embed $F_q \rightarrow F_{q^2}$, so without loss of generality we can work over $F_{q^2}$. Then define $\alpha: x \rightarrow x^q$ to be an order two involution, the $r^{th}$ power of Frobenius if $q=p^r$. Then an operator is unitary if the inner product $\langle x, y\rangle = \sum_i x_i^q y_i$ on an $F_{q^2}$ vector space is preserved. In other words, if $UU^* = I_N$ where $U^* = (U^\alpha)^T$. \\

When $p \nmid N$, we have an $N^{th}$ root of unity $\alpha$, so the unnormalized DFT is $A_{ij} = \alpha^{ij}$. $\sqrt{N}$ may not be in $F_{q^2}$. If not, form the quadratic extension $F_{q^2}[x]/(x^2-N) \cong F_{q^4}$ which surely contains $\sqrt{N}$.

$$\left(AA^*\right)_{ik} = \sum_{j=0}^{n-1} \alpha^{j(i+qk)} = N\delta_{i,-qk}$$

where we use a geometric sum to simplify. Normalizing, $UU^* = \delta_{i,-qk}$. Thus $U$ is unitary iff $q \equiv -1 (\mod N)$. \\

\section{Unitarity for the symmetric group}

\subsection{Weyl's unitary trick} The Fourier transform described in \S \ref{fourier-transform-nonabelian-groups} is not unitary. In order make the transform unitary, we need to include factors involving square roots. We'll start by working over the complex numbers, and then consider more economical subfields. Consider

$$\hat{f}(\tau) = \sum_{g \in G} \sqrt{\frac{d_\tau}{|G|}} f(g)\tau(g)$$

where every $\tau(g) \in U(n)$ is now a unitary matrix, and $d_\tau$ is the dimension of the representation $\tau$. The overall DFT is now unitary with these two modifications \cite{beals-qft-symmetric-group}. Note that every representation of a finite group, or even a compact Lie group, may be made unitary by Weyl's unitary trick. Recall the invariant inner product $\langle \cdot , \cdot \rangle_G = \frac{1}{|G|}\sum_{g \in G} \langle \rho(g)x, \rho(g)y \rangle$. By using the change-of-basis matrix associated to this inner product, we can conjugate $\rho$ to a unitary representation $\tau$. There are different versions of the unitary trick, but we will use the following. Let

$$P = \int_G \rho(g)\rho(g)^* dg$$

and $P=Q^2$, $Q$ is a principal square root of $P$, and $dg$ is the Haar measure. $Q$ is the change-of-basis matrix required, and $\tau(g) = Q^{-1}\rho(g)Q$ is unitary. For finite groups, $P = \frac{1}{|G|}\sum_{g \in G} \rho(g)\rho(g)^*$. $Q$ can be computed by diagonalizing the positive symmetric matrix $P$. Let $P = MDM^{-1}$, then $Q = MD^{1/2}M^{-1}$.

\subsection{Over number fields}

For the symmetric group, the characteristic-zero case is immediate from Young's orthogonal form. Its matrices are real and orthogonal, and hence unitary over $\mathbb{C}$. Their entries, together with the normalization factors $\sqrt{d_\tau/|G|}$ appearing in the Fourier transform, lie in a finite extension of $\mathbb{Q}$. Thus the unitary DFT is defined over a number field containing these square roots.

The situation over finite fields requires a different argument, since there is no notion of positivity and an invariant Hermitian form need not admit an ordinary Cholesky decomposition.

\subsection{Unitary representations over finite fields}

Let $F=GF(q^2)$ with its usual conjugation $x\mapsto x^q$, and assume that the characteristic of $F$ does not divide $n!$. We find a change of basis which makes each Specht representation unitary in two steps: first compute the Gram matrix $U$ of an $S_n$-invariant Hermitian form, and then factor $U=A^\ast A$ by an extended Cholesky decomposition. \\

A $G$-invariant Hermitian form satisfies $\langle gx, gy\rangle = \langle x, y \rangle$ for all $g \in G$. Using the convention $\langle x,y\rangle=x^\ast Uy$, we obtain

$$\rho(g)^\ast U\rho(g)=U.$$

We regard the $d_\rho^2$ entries $u_{ij}$ of $U$ as variables over $F$ and stack the resulting linear equations for the elements $g\in G$. Since the Specht representations are defined over the prime subfield, conjugation fixes the entries of $\rho(g)$. In the semisimple case the space of invariant forms is one dimensional, and we select a nonsingular Hermitian solution.

\begin{lstlisting}
from sage.rings.polynomial.polynomial_ring_constructor import PolynomialRing
G = Permutations(self._n)
F = self._ring
rho = self._specht.representation_matrix

# compute the invariant Hermitian matrix
d_rho = self._specht.dimension()
R = PolynomialRing(F, 'u', d_rho**2)
U_vars = R.gens()
Utemp = matrix(R, d_rho, d_rho, U_vars)

def augmented_matrix(g):
      rho_g = rho(g)
      equation_matrix = rho_g.conjugate_transpose() * Utemp * rho_g - Utemp
      augmented_system = []
      for i in range(d_rho):
           for j in range(d_rho):
               linear_expression = equation_matrix[i, j]
               row = [linear_expression.coefficient(u) for u in U_vars]
               augmented_system.append(row)
        return augmented_system

total_system = sum((augmented_matrix(g) for g in G), [])
null_space = matrix(F, total_system).right_kernel()
U = matrix(F, d_rho, d_rho, null_space.basis()[0])
A = U.cholesky(extended=True).conjugate_transpose()
return A
\end{lstlisting}

The extended Cholesky decomposition over a finite field uses Gaussian elimination adapted to finite-field conjugation. The implementation was adapted from GAP's \texttt{BaseChangeToCanonical} method for Hermitian forms. \\

Then the representations $\tilde{\rho}(g)=A\rho(g)A^{-1}$ are unitary since

$$\tilde{\rho}(g)^\ast\tilde{\rho}(g)
= A^{-\ast}\rho(g)^\ast A^\ast A\rho(g)A^{-1}$$
$$=A^{-\ast}\rho(g)^\ast U\rho(g)A^{-1}
=A^{-\ast}UA^{-1}$$
$$=A^{-\ast}A^\ast AA^{-1}=I.$$

\end{document}